\documentclass[11pt]{article}

\usepackage{fullpage}
\usepackage{amsmath, graphicx, amsthm, amsfonts, amssymb, amstext, mathrsfs}
\usepackage{ragged2e, enumerate, graphicx, lscape, framed}
\usepackage{pgf, tikz, float, subcaption}
\usetikzlibrary{graphs, graphs.standard, decorations.pathreplacing, calligraphy}
\usepackage[symbols,nogroupskip,sort=none]{glossaries-extra}
\usepackage{subfiles}
\usepackage[T1]{fontenc}

\newtheorem{theorem}{Theorem}[section]
\newtheorem{lemma}[theorem]{Lemma}
\newtheorem{proposition}[theorem]{Proposition}
\newtheorem{conjecture}[theorem]{Conjecture}
\newtheorem{corollary}[theorem]{Corollary}

\usepackage{hyperref}
\hypersetup{
	colorlinks=true,
    urlcolor=purple,
	linkcolor=blue,
    citecolor=purple,
}

\DeclareMathOperator{\supp}{supp}
\DeclareMathOperator{\ex}{ex}
\DeclareMathOperator{\mex}{mex}

\title{Localization of spectral Tur\'{a}n-type theorems}
\author{M. Rajesh Kannan, Hitesh Kumar, Shivaramakrishna Pragada}
\date{}

\begin{document}
\maketitle

\begin{abstract} Let $G$ be a graph, and let $v$ and $e$ be a vertex and an edge of $G$, respectively. Define $c(v)$ (resp. $c(e)$) to be the order of the largest clique in $G$ containing $v$ (resp. $e$). Denote the adjacency eigenvalues of $G$ by $\lambda_1 \ge \cdots \ge \lambda_n$. We study localized refinements of spectral Tur\'{a}n-type theorems by replacing global parameters such as the clique number $\omega(G)$, size $m$ and order $n$ of $G$ with local quantities $c(v)$ and $c(e)$.

Motivated by a conjecture of Elphick, Linz and Wocjan (2024), we first propose a vertex-localized strengthening of Wilf's inequality:
\[ \sqrt{s^{+}(G)} \le \sum_{v\in V(G)}\left(1-\frac{1}{c(v)}\right), \]
where $s^+(G) = \sum_{\lambda_i > 0}\lambda_i^2$. 
Inspired by the Bollob\'{a}s-Nikiforov conjecture (2007) on the first two eigenvalues, we then introduce an edge-localized analogue:
\[\lambda_1^2(G) + \lambda_2^2(G) \le \sum_{e\in E(G)} 2\left(1-\frac{1}{c(e)}\right).\]
As evidence of their validity, we verify the above conjectures for diamond-free graphs and random graphs. We also propose strengthening of the spectral versions of the Erd\H{o}s-Stone-Simonovits Theorem by replacing the spectral radius with $\sqrt{s^{+}(G)}$ and establish it for all $F$-free graphs with $\chi(F)=3$. A key ingredient in our proofs is a general upper bound relating $\sqrt{s^{+}(G)}$ to the triangle count $t(G)$. Finally, we prove a localized version of Nikiforov's walk inequality and conjecture a stronger localized version. These results contribute to the broader program of localizing spectral extremal inequalities.
\end{abstract}



\section{Introduction}

\subsection{Notation}

We use standard graph theory notation and terminology throughout the paper. We assume all graphs are simple and connected unless specified otherwise. Let $G = (V(G), E(G))$ be a graph of order $n$ and size $m$. Let $t(G)$, $\Delta(G)$, $\omega(G)$ and $\chi(G)$ denote the number of triangles, maximum degree, clique number and chromatic number of $G$, respectively. For a graph $F$, we say that $G$ is $F$-free if $G$ does not contain $F$ as a subgraph (not necessarily induced). Let $\deg(v)$ denote the degree of vertex $v$ in $G$. For an edge $e\in E(G)$ (or a vertex $v\in V(G)$), we will denote by $c(e)$ (resp. $c(v)$) the order of the largest clique in $G$ containing the edge $e$ (resp. the vertex $v$). The \textit{adjacency matrix} of $G$ is an $n\times n$ matrix $A(G) = [a_{uv}]$, where $a_{uv} = 1$ if $u$ and $v$ are adjacent, and $0$ otherwise. The \emph{eigenvalues} of $G$ are the eigenvalues of $A(G)$. Since $A(G)$ is a real symmetric matrix, all eigenvalues of $A(G)$ are real and can be listed as $\lambda_1(G) \geq  \cdots \geq \lambda_n(G)$. Let  
\[s^+(G)=\sum_{\lambda_i>0}\lambda_i^2(G) \text{ and }s^-(G)=\sum_{\lambda_i<0}\lambda_i^2(G).\]
We call $s^+(G)$ (resp. $s^-(G)$) the \emph{positive} (resp. \emph{negative}) \emph{square energy} of $G$. Let $n^+$ and $n^-$ denote the number of positive and negative eigenvalues of $A(G)$, respectively.

For a graph $G$, we can assign a (non-negative) \emph{weight function} $w: E(G)\rightarrow \mathbb{R}^{\ge 0}$ to its edges and define a \emph{weighted adjacency matrix} $W(G) = [w_{uv}]$, where $w_{uv} = w(uv)$ if $u$ and $v$ are adjacent, and $w_{uv} = 0$ otherwise. For an integer $r\ge 1$ and a vertex $v\in V(G)$, denote by $w_r(v)$ the number of walks with $r$ vertices starting from $v$ in $G$. The total number of $r$-walks in $G$ is then given by 
\[w_r(G) = \sum_{v\in V(G)} w_r(v).\] 

\subsection{Known results}
One of the central themes in graph theory is understanding how local constraints influence global structure. A classic example is determining the maximum number of edges a graph can have while avoiding a specific forbidden subgraph. This area of study, known as extremal graph theory, seeks to characterize such optimal graphs under given restrictions. A foundational result in this field is Tu\'{r}an's Theorem, which provides an exact bound on the number of edges in a graph that does not contain a complete subgraph $K_{\omega+1}$ .

\begin{theorem}[Tu\'{r}an's Theorem]\label{thm:Turan}
For any graph $G$ with the clique number $\omega$, we have
\begin{equation*}
   m \le \frac{n^2}{2}\left( 1 - \frac{1}{\omega}\right). 
\end{equation*}
Equality holds if and only if $G$ is a complete regular $\omega$-partite graph.
\end{theorem}

Wilf \cite{wilf1986spectral} obtained the following spectral generalization of Tu\'{r}an's result. 

\begin{theorem}[Wilf's Inequality \cite{wilf1986spectral}]\label{thm:Wilf}
For any graph $G$, we have
\begin{equation*}
    \lambda_1(G) \le n \left( 1 - \frac{1}{\omega}\right).
\end{equation*}   
Equality holds if and only if $G$ is a complete regular $\omega$-partite graph.
\end{theorem}

 Wilf's Inequality is stronger than Tu\'{r}an's Theorem since the spectral radius is always larger than the average degree of the graph, i.e.,  $\lambda_1 \ge \frac{2m}{n}$.

Furthermore, Edwards and Elphick \cite{Edwards_Elphick_1983} conjectured the following generalization, which was later proved by Nikiforov \cite{Nikiforov_2002}, now known as the spectral Tu\'{r}an's Theorem. 

\begin{theorem}[Spectral Tu\'{r}an's Theorem \cite{Nikiforov_2002}]\label{thm:spectral_Turan}
For any graph $G$, we have
\begin{equation*}
    \lambda_1^2(G) \le 2m \left( 1 - \frac{1}{\omega}\right).
\end{equation*}
\end{theorem}

The following extension of Tur\'{a}n's Theorem to $F$-free graphs is well-known.

\begin{theorem}[Erd\H{o}s-Stone-Simonovits Theorem \cite{Erdos_Simonovits_1966, Erdos_Stone_1946}]
For any $F$-free graph $G$ where $\chi(F)\ge 3$, we have 
\[ m\le \frac{n^2}{2}\left( 1 - \frac{1}{\chi(F)-1} + o(1)\right).\]    
\end{theorem}

Nikiforov \cite{Nikiforov_Erdos_stone_2009} established the spectral version of the above theorem in the spirit of Wilf's Inequality. 

\begin{theorem}[\cite{Nikiforov_Erdos_stone_2009}]\label{thm:spectral_ESS_Vertex}
    For any $F$-free graph $G$ where $\chi(F) \ge 3$, we have
    \[\lambda_1(G)\le n \left(1 - \frac{1}{\chi(F)-1} + o(1)\right).\]
\end{theorem}

Li, Liu and Zhang \cite{li_liu_zhang_2025_stone_simonovits} proved an edge-spectral generalization of Erd\H{o}s-Stone-Simonovits Theorem.

\begin{theorem}[\cite{li_liu_zhang_2025_stone_simonovits}]\label{thm:spectral_ESS_Edge}
    For any $F$-free graph $G$ where $\chi(F) \ge 3$, we have
    \[\lambda_1^2(G) \le 2m\left(1 - \frac{1}{\chi(F)-1} + o(1)\right).\]
\end{theorem}

Recently, there has been an interest in generalizing classical results by assigning a weight function to the edges (resp. vertices) of the graph and then proving a bound on this weight function which implies the classical result. This is called \emph{edge localization} (resp. \emph{vertex localization}). Malec and Tompkins \cite{Malec_Tompkins_2023} edge-localized a wide range of classical results in graph theory, and many interesting follow-up papers have appeared, see \cite{Adak_Chandran_2025_erdos_gallai, Balogh_Bradac_Lidicky_2025, Kirsch_Nir_2024, Zhao_Zhang_2025_Erdos_Gallai, Zhao_Zhang_2025_long_cycles}. Tu\'{r}an's Theorem was edge-localized independently by Brada\v{c} \cite{Bradac_2022} and Malec and Tompkins \cite{Malec_Tompkins_2023}. Recently, Adak and Chandran \cite{Adak_Chandran_2025_Turan} vertex-localized Tu\'{r}an's Theorem. Liu and Ning \cite{Liu_Ning_2026} obtained the following edge localization of the spectral Tu\'{r}an's Theorem. 

\begin{theorem}[Edge-localized Spectral Tu\'{r}an's Theorem \cite{Liu_Ning_2026}] \label{thm:spectral_edge_localization_Turan}
For any graph $G$, we have
\begin{equation*}
   \lambda_1^2(G) \le \sum_{e\in E(G)} 2\left(1-\frac{1}{c(e)}\right). 
\end{equation*}
\end{theorem}

We note here that in Theorems \ref{thm:spectral_Turan} and \ref{thm:spectral_edge_localization_Turan}, equality holds if and only if $G$ is a complete bipartite graph when $\omega=2$ or a regular complete $\omega$-partite graph when $\omega\ge 3$. 

Recently, Liu and Ning generalized Theorem \ref{thm:spectral_edge_localization_Turan} to weighted graphs, which implies the results of Brada\v{c}, and Malec and Tompkins for Tu\'{r}an's Theorem.

\begin{theorem}[Weighted Edge-localized Spectral Tu\'{r}an's Theorem \cite{Liu_Ning_2025_weighted}] \label{thm:weighted_spectral_edge_localization_Turan}
For any connected graph $G$ with weighted adjacency matrix $W(G)$, we have
\begin{equation*}
   \lambda_1^2(W(G)) \le \sum_{e\in E(G)} 2\left(1-\frac{1}{c(e)}\right)w(e)^2. 
\end{equation*}
\end{theorem}

\subsection{Our contribution}

Elphick, Linz and Wocjan \cite{Elphick_Linz_Wocjan_2024} proposed the following generalization of Wilf's Inequality. 

\begin{conjecture}[\cite{Elphick_Linz_Wocjan_2024}]\label{conj:splus_omega}
For any graph $G$, we have 
    \[\sqrt{s^+(G)} \leq n\left(1 - \frac{1}{\omega} \right).\]
\end{conjecture}

In \cite{Elphick_Linz_Wocjan_2024}, the authors have verified the above conjecture for weakly-perfect graphs, Kneser graphs, triangle-free graphs and random graphs. We prove a weaker but non-trivial upper bound for $\sqrt{s^+}$ of any graph $G$ as follows.

\begin{theorem}\label{thm:splus_weak_2omegaminus2}
Let $G$ be a graph of order $n$ and the clique number $\omega$. Then
\[\sqrt{s^+(G)} \leq n\sqrt{1-\frac{1}{\omega}-\frac{1}{\omega^2}} \ .\] 
\end{theorem}

Next, we propose a vertex-localized version of the above conjecture.

\begin{conjecture}\label{conj:vertex_local_wilf}
For any graph $G$, we have
    \[\sqrt{s^+(G)}\le \sum_{v\in V} \left(1-\frac{1}{c(v)}\right).\]
\end{conjecture}

We present a case for this conjecture by proving it for diamond-free graphs (Theorem \ref{thm:local_wilf_diamond_free}), and almost all graphs using Erd\H{o}s-R\'enyi random graph $G(n,1/2)$ (Theorem \ref{thm:local_wilf_random}). We have also computationally verified this conjecture for graphs up to 9 vertices \href{https://github.com/Shivaramkratos/Code-for-Localization/blob/main/Local_wilf_check.sage}{here}. To prove Conjecture \ref{conj:vertex_local_wilf} for diamond-free graphs, we first establish the following general upper bound for $s^+$ in terms of the number of triangles $t(G)$, which we believe is of independent interest.

\begin{theorem}\label{thm:splus_triangles}
    Let $G$ be a graph of order $n$, size $m$, and $t(G)$ be the number of triangles in $G$. Then
    \[\sqrt{s^+(G)} \leq \frac{n}{2} + \frac{3t(G)}{\lambda^2_1(G)}.\]
    Equality holds if and only if $G$ is a complete bipartite graph.
\end{theorem}

Contrast this with a known result of Elphick and Linz \cite{Elphick_Linz_2024_symmetry} which states that $\sqrt{s^-(G)}\le \frac{n}{2}$ for any graph $G$ of order $n$. This further suggests an asymmetry between $s^+$ and $s^-$; also refer \cite{Zhang_2024_square_energy}. 

For the $F$-free graphs, we propose a strengthening of Nikiforov's Theorem \ref{thm:spectral_ESS_Vertex} by replacing $\lambda_1$ with $\sqrt{s^+}$. 

\begin{conjecture}\label{conj:Erdos_Stone_wilf}
For any $F$-free graph $G$ where $\chi(F) \ge 3$, we have
    \[\sqrt{s^+(G)}\le n \left(1 - \frac{1}{\chi(F)-1} + o(1)\right).\]
\end{conjecture}

We prove this conjecture for $F$-free graphs when $\chi(F)=3$ (Theorem \ref{thm:ESS_general_F_free_chi_3}).

Bollob\'{a}s and Nikiforov \cite{bollobas2007cliques} proposed to improve Theorem \ref{thm:spectral_Turan} by introducing $\lambda_2$.

\begin{conjecture}[Bollob\'{a}s - Nikiforov Conjecture \cite{bollobas2007cliques}]\label{conj:BN}
For any graph $G \neq K_n$, we have 
\[\lambda_1^2(G) + \lambda^2_2(G) \leq 2m\left(1 - \frac{1}{\omega} \right).\]
\end{conjecture}

The above conjecture has been verified for weakly perfect graphs by Ando and Lin \cite{ando_lin_2015}, triangle-free graphs by Lin, Ning and Wu \cite{Lin_Ning_Wu_2021}, regular graphs by Zhang \cite{zhang2024regular}, and random graphs by Liu and Bu \cite{Liu_Bu_2025}. Kumar and Pragada \cite{Bollobas_Nikiforov_few_triangles_2025} verified the conjecture for graphs with few triangles, which includes cycle-free graphs, book-free graphs and planar graphs. Zeng and Zhang \cite{Zeng_Zhang_2025} proved the conjecture for line graphs and also strengthened some results from \cite{Bollobas_Nikiforov_few_triangles_2025}. Motivated by the recent localization results, it is only natural to ask if the following edge-localized refinement of the Bollob\'{a}s-Nikiforov Conjecture is true.

\begin{conjecture}\label{conj:local_BN}
For any graph $G \neq K_n$, we have
    \[\lambda_1^2(G) + \lambda_2^2(G) \le \sum_{e\in E(G)} 2\left(1-\frac{1}{c(e)}\right).\]
\end{conjecture}

Clearly, Conjecture \ref{conj:local_BN} is true for triangle-free graphs. As evidence for this conjecture, we prove it for the next non-trivial case of diamond-free graphs (Theorem \ref{thm:local_BN_diamond_free}) whenever $t(G)\notin \{1,2,3,4\}$, and random graphs (Theorem \ref{thm:local_BN_random}). We have computationally verified this conjecture for graphs up to 9 vertices \href{https://github.com/Shivaramkratos/Code-for-Localization/blob/main/BN_check.sage}{here}.  

We note that equality in Conjecture~\ref{conj:local_BN} is attained not just by regular complete multipartite graphs (barring $K_n$), but also by the disjoint union of two regular complete multipartite graphs of not necessarily the same regularities. This suggests that Conjecture~\ref{conj:local_BN} captures a rich extremal phenomenon.


Recently, Li, Liu and Zhang \cite{li_liu_zhang_2025_color_critical} proposed a strengthening of the Bollob\'{a}s-Nikiforov Conjecture for $F$-free graphs where $F$ is a color-critical graph.

\begin{conjecture}[\cite{li_liu_zhang_2025_color_critical}] Let $F$ be a color-critical graph with $\chi(F)\ge 4$. Then any $F$-free graph $G$ with sufficiently large size $m$ satisfies
\[\lambda_1^2(G) + \lambda^2_2(G) \leq 2m\left(1 - \frac{1}{\chi(F)-1}\right).\]
\end{conjecture}

We propose a weaker bound, but for any graph $F$ with $\chi(F)\ge 3$.

\begin{conjecture}
 For any $F$-free graph $G$ where $\chi(F) \ge 3$, we have 
\[\lambda_1^2(G) + \lambda^2_2(G) \leq 2m\left(1 - \frac{1}{\chi(F)-1} + o(1) \right).\]
\end{conjecture}

We prove this statement when $\chi(F) = 3$ (Theorem \ref{thm:BN_F_free_chi_3}). 
 
Nikiforov \cite{Nikiforov_2002} gave an upper bound on $\lambda_1$ using the number of $r$-walks in the graph. 

\begin{theorem}[Nikiforov's Walk Inequality \cite{Nikiforov_2002}]\label{thm:Nikiforov_walk} Let $G$ be a graph. For any integer $r\ge 1$, we have 
\begin{equation*}
\lambda_1^{r}(G) \le w_{r}(G) \left( 1 - \frac{1}{\omega}\right).
\end{equation*}
Equality holds if and only if one of the following holds:
\begin{enumerate}[$(i)$]
    \item If $\omega >2$, then $G$ is a regular $\omega$-partite graph.  
    \item If $\omega =  2$, then $G$ is a complete bipartite graph, and if $r$ is odd, then $G$ is regular. 
\end{enumerate}
\end{theorem}

When $r= 1$, Theorem \ref{thm:Nikiforov_walk} implies Theorem \ref{thm:Wilf} since $w_1(G) = n$. When $r = 2$, Theorem \ref{thm:Nikiforov_walk} implies Theorem \ref{thm:spectral_Turan} since $w_2(G) = 2m$. We obtain a localized strengthening of Nikiforov’s inequality.

\begin{theorem}\label{thm:partial_walks}
Let $G$ be a graphFor any integer $r\ge 1$, we have 
\[ \lambda_1^r(G) \le \sum_{v\in V(G)} w_r(v)\sqrt{1-\frac{1}{c(v)}}\sqrt{1-\frac{1}{\omega (G)}}.\]
The equality case is similar to that of Theorem \ref{thm:Nikiforov_walk}.
\end{theorem}

Although Theorem \ref{thm:partial_walks} is localized, we believe that the following stronger localized version is true. 

\begin{conjecture}\label{conj:local_walks}
For any integer $r\ge 1$, we have 
\[ \lambda_1^r(G) \le \sum_{v\in V(G)} w_r(v)\frac{c(v)-1}{c(v)}.\]
\end{conjecture}

\subsection{Organization of the paper}

This paper is organized as follows. In Section~\ref{sec:Preliminaries}, we collect some known results about Motzkin-Straus type theorems, majorization, triangle-counting and random graphs, that will be useful later. Section~\ref{sec:Wilf} is devoted to Wilf-type results while in Section~\ref{sec:BN} we deal with generalizations of Bollob\'{a}s-Nikiforov conjecture. Finally, in Section~\ref{sec:Walk_inequalities} we focus on inequalities involving walks.

\section{Preliminaries}
\label{sec:Preliminaries}

\subsection{Weighted Motzkin-Straus type theorems}

Consider a graph $G = (V(G), E(G))$ with $\vert V \vert$, the number of vertices, equals $n$. Let $S$ denote the standard simplex given by 
\[S = \left\{x \in \mathbb{R}^n : \sum_{v \in V(G)} x_v = 1, x_v \geq 0, v\in V(G)\right\}.\]
For a vector $x\in \mathbb{R}^n$, the \emph{support} of $x$ is given by $\supp(x) = \{v\in V(G) : x_v\neq 0\}$. Motzkin-Strauss \cite{Motzkin_Straus_1965} established the following remarkable result. 

\begin{theorem}[\cite{Motzkin_Straus_1965}]\label{thm:Motzkin_Strauss}
	Let $G$  be a graph with adjacency matrix $A(G)$ and clique number $\omega$. For any $x \in S$, we have  
	\[x^TA(G)x \leq 1- \frac{1}{\omega}.\]
    Equality holds if and only if $\supp(x)$ induces a complete $\omega$-partite graph.
\end{theorem}

Instead of working with $A(G)$, one can also work with a weighted (non-negative weights) adjacency matrix $W(G) = [w_{uv}] \in \mathbb{R}^{n \times n}$ of the given graph $G$, and maximize the quantity $F_G(x):= x^TW(G)x$ on the simplex $S$. For the weight matrix $W(G)=[w_{uv}]$, where $w_{uv} = \frac{c(uv)}{c(uv) - 1}$, it was shown in \cite{Bradac_2022, Liu_Ning_2026}) that $F_G(x)\le 1$ for all $x\in S$ and equality holds if and only if $\supp(x)$ induces a complete $\omega$-partite graph. Recently, Liu and Ning \cite{Liu_Ning_2025_weighted} established the following.

\begin{lemma}[\cite{Liu_Ning_2025_weighted}]\label{lemma:weighted_M_S}
    Let $W(G)=[w_{uv}]$ be a weighted matrix where $w_{uv} = \frac{1}{2}\left(\frac{c(u)}{c(u)-1} + \frac{c(v)}{c(v)-1}\right)$. Then $F_G(x)\le 1$ for all $x\in S$, and equality holds if and only if $\supp(x)$ induces a complete $\omega$-partite graph whose vertex classes $V_1, \ldots, V_{\omega}$ satisfy $\sum_{v\in V_i}x_v = \frac{1}{\omega}$. 
\end{lemma}

We observe the following corollary, which we require later.

\begin{lemma}\label{lemma:weighted_M_S_2}
     Let $W(G)=[w_{uv}]$ be a weighted matrix where $w_{uv} = \sqrt{\frac{c(u)c(v)}{(c(u) - 1)(c(v)-1)}}$. Then $F_G(x)\le 1$ for all $x\in S$, and equality holds if and only if $\supp(x)$ induces a complete $\omega$-partite graph whose vertex classes $V_1, \ldots, V_{\omega}$ satisfy $\sum_{v\in V_i}x_v = \frac{1}{\omega}$. 
\end{lemma}

\begin{proof} Let $x\in S$. Using AM-GM inequality and Lemma \ref{lemma:weighted_M_S}, we get
\[ x^T W(G)x = \sum_{uv \in E(G)} 2\sqrt{\frac{c(u)c(v)}{(c(u) - 1)(c(v)-1)}} x_ux_v \le \sum_{uv\in E(G)} \left(\frac{c(u)}{c(u)-1} + \frac{c(v)}{c(v)-1}\right)x_ux_v \le 1.\]
It is easily seen that the equality case is similar to that of Lemma \ref{lemma:weighted_M_S}.
\end{proof}

\subsection{Majorization}\label{subsec:majorization}

For a vector $x \in \mathbb{R}^n$, let $x^\downarrow$ denote the vector obtained by rearranging the entries of $x$ in the non-increasing order. Given two vectors $x, y \in \mathbb{R}^n $, we say that  $x$ is \textit{weakly majorized} by $y$, if
\begin{gather*}
    \sum_{j=1}^{k} {x_j}^\downarrow \leq \sum_{j=1}^{k} {y_j}^\downarrow
\end{gather*}
for all $1 \leq k \leq n$, and is denoted by $x \prec_w y$. 
Moreover, if $x \prec_w y \ \text{and}  \ \sum_{i=1}^{n} x_i^\downarrow = \sum_{i=1}^{n} y_i^\downarrow,$
then we say that $x$ is \textit{majorized} by $y$, and is denoted by $x \prec y$. 
For $x = (x_1,x_2,\dots,x_n) \in \mathbb{R}^n$ and $1 \leq p < \infty$, define $\Vert x \Vert_p = \Big\{\sum_{i=1}^{n} \vert x_i \vert^p \Big\}^\frac{1}{p}.$

\begin{theorem}[{\cite{Lin_Ning_Wu_2021}}]\label{thm:majorization}
    Let $x = (x_1,\dots,x_n),\ y = (y_1,\dots,y_n) \in \mathbb{R}^n_{\geq 0}$. If $y \prec_w x$, then 
    \[\Vert y\Vert_p \leq \Vert x\Vert_p\]
    for any real number $p> 1$, and equality holds if and only if $x=y$.
\end{theorem}

\subsection{Triangle-counting}

One way to count triangles in a graph is by counting the number of edges in vertex-neighbourhoods.

\begin{lemma}\label{lemma:triangle_counting} For every graph $G$, the number of triangles $t(G)$ in $G$ is given by
\[t(G) =  \frac{1}{3}\sum_{v\in V(G)} m(G[N(v)]),\]
where $m(G[N(v)])$ denotes the number of edges in the subgraph induced by $N(v)$.
\end{lemma}

Let $F$ and $H$ be graphs. Let $\ex(n, H, F)$ denote the maximum number of copies of the graph $H$ in an $F$-free graph of order $n$. We mention the following general result concerning $F$-free graphs. 

\begin{theorem}[\cite{gerbner2019Ffreecounting}, Theorem 2]\label{thm:n_general_counts}
    Let $H$ and $F$ be graphs and assume $\chi(F) = k$. Then
    \[\ex(n,H,F) \le \ex(n,H,K_k) + o(n^{|V(H)|}).\]
\end{theorem}

Let $\mex(m, H, F)$ denote the maximum number of copies of the graph $H$ in an $F$-free graph of size $m$. We state a similar counting result as above, but in terms of the size.  

\begin{theorem}[\cite{wang2025generalizedturan}, Corollary 1.2]\label{thm:m_general_counts}
Let $r \ge 3$ be an integer and $F$ be a graph. If $\ex(n, K_r, F ) = o(n^s)$ for some $1 < s \le r$, then $\mex(m, K_r, F) = o\left(m^{\frac{(r-1)s}{r+s-2}}\right)$.
\end{theorem}

We require the following corollary of the above theorems when $H$ is a triangle. 

\begin{corollary}\label{cor:triangles_F_free}
Let $F$ be a fixed graph with $\chi(F)=3$. Then 
\begin{enumerate}[$(i)$]
    \item $\ex(n, K_3, F) \le o(n^{3})$
    \item $\mex(m, K_3, F) \le o(m^{3/2})$
\end{enumerate}
\end{corollary}

\begin{proof}
For $(i)$, taking $H = K_3$ in Theorem \ref{thm:n_general_counts} and noting that $\ex(n,K_3,K_3) = 0$, completes the proof. For $(ii)$, by the first assertion, we have $\ex(n,K_3,F) \le o(n^3)$, implying that for some $\epsilon>0$, we have $\ex(n,K_3,F) = o(n^{3-\epsilon})$. Now taking $r=3$ and $s = 3-\epsilon$ in Theorem \ref{thm:m_general_counts}, gives $\mex(m,K_3,F) = o\left(m^{\frac{6-2\epsilon}{4-\epsilon}}\right) \le o(m^{3/2})$.
\end{proof}

We also require the following well-known result of Erd\H{o}s and Gallai \cite{Erdos_Gallai_1959}.

\begin{theorem}[Erd\H{o}s-Gallai \cite{Erdos_Gallai_1959}]\label{thm:erdos_gallai} For a graph $G$ of order $n$ and size $m$ with no $P_k$ as a subgraph, we have
\[m\le \frac{n(k-2)}{2}.\]
\end{theorem}

\subsection{Random graphs}

For $n\in \mathbb{N}$ and $0<p<1$, the Erd\"{o}s-Renyi random graph $G(n,p)$ is a probability space over the set of graphs with vertex set $V(G) = \{v_1, \ldots, v_n \}$ such that whenever $i<j$, the event that $v_i$ is adjacent to $v_j$ occurs with probability $p$, and these events are mutually independent. 

Consider a symmetric matrix $A = [a_{ij}]$ of order $n$ such that $a_{ii} = 0$ for all $i$ and $a_{ij}$'s are independent random variables whenever $i < j$ with probability mass function given by: 
\[\mathbb{P}(a_{ij} = 1) = p \text{ and } \mathbb{P}(a_{ij} = 0) = 1-p,\]
where $0 < p <1$. Then $A$ is the adjacency matrix associated with the random graph $G(n, p)$. The eigenvalues of $G(n,p)$ are then defined to be the eigenvalues of $A$. 

We say that a property $\Lambda$ holds for $G(n,p)$ \emph{aymptotically almost surely} (a.a.s. for short) if 
\[\mathbb{P}(\Lambda \text{ holds for }G(n,p)) \rightarrow 1 \text{ as } n\rightarrow \infty.\] 

The following is known about the first and the second eigenvalue of random graphs. 

\begin{theorem}\label{thm:random_graph_eigenvalues} For the random graph $G = G(n,p)$, the following hold a.a.s: 
\begin{enumerate}[$(i)$]
    \item \emph{(\cite{Juhasz_1981})}  $\lambda_1 (G) = (p + o(1))n$.
    \item \emph{(\cite{Furedi_Komlos_1981})}  $\lambda_2(G) \le 2\left(\sqrt{p(1-p)} + o(1)\right)\sqrt{n}$.
\end{enumerate}
\end{theorem}

The clique number of random graphs is given by the following well-known result.

\begin{theorem}[\cite{Grimmett_McDiarmid_1975, Bollobas_Erdos_1976}] \label{thm:random_graphs_clique_no} For the random graph $G = G(n,p)$, a.a.s we have 
\[\omega(G) = \left( \frac{2}{\log(\frac{1}{p})} + o(1)\right)\log n.\]
\end{theorem}

Using Wigner's semicircle law \cite{Wigner_1958}, one can determine $s^+$ and $s^-$ of the random graph $G(n, 1/2)$ as was done in \cite{Elphick_Linz_2024_symmetry}.  

\begin{theorem}[\cite{Elphick_Linz_2024_symmetry}]\label{thm:random_graph_square_energy}
For the random graph $G = G(n,p)$, a.a.s we have 
\[ s^+(G) = \frac{1}{2}\left( \frac{3}{4} + o(1) \right)n^2 \text{ and } s^-(G) = \frac{1}{2}\left(\frac{1}{4} + o(1)\right)n^2.\]
\end{theorem}

\section{Wilf-type inequalities}
\label{sec:Wilf}

\subsection{Some general upper bounds}

Firstly, we give an upper bound on $\sqrt{s^+}$ in terms of the number of triangles in the graph, thus proving Theorem \ref{thm:splus_triangles}.

\begin{proof}[Proof of Theorem \ref{thm:splus_triangles}]
It is well-known that
\[ \sum_{i=1}^n \lambda_i^3 = \sum_{\lambda_i>0} \lambda_i^3 - \sum_{\lambda_i<0} |\lambda_i|^3 =  6t(G).\]
We have
\[ s^- \ge  \frac{\sum_{\lambda_i<0} |\lambda_i|^3}{|\lambda_n|} = \frac{\sum_{i=1}^{n^+} \lambda_i^3 -6t(G)}{|\lambda_n|} \geq \frac{\sum_{i=1}^{n^+} \lambda_i^3 -6t(G)}{\lambda_1} \ge \lambda_1^2 - \frac{6t(G)}{\lambda_1}, \]
which gives
\begin{equation}\label{eq:triangle_s_minus}
\lambda_1 \leq \sqrt{s^- + \frac{6t(G)}{\lambda_1}}.    
\end{equation}
Thus,
\begin{equation*}
    \sqrt{s^+} = \frac{\lambda_1}{\lambda_1}\sqrt{s^+} \le \frac{1}{\lambda_1}\sqrt{s^- +\frac{6t(G)}{\lambda_1}}\sqrt{s^+}\le \frac{1}{\lambda_1}\frac{2m +\frac{6t(G)}{\lambda_1}}{2} \le \frac{n}{2} + \frac{3t(G)}{\lambda^2_1},
\end{equation*}
proving the desired inequality.

It is clear from above that if equality holds, then $\lambda_1 = |\lambda_n|$, implying $G$ is bipartite. Thus $s^+ = s^- = m$. Now, $\sqrt{s^+} = \frac{n}{2}$ only if $G$ is complete bipartite. This completes the proof. 
\end{proof}

We note here that inequality \eqref{eq:triangle_s_minus} in the above proof can be rearranged to obtain a lower bound on the number of triangles in a graph:
\begin{equation}\label{eq:triangle_s_minus_rearranged}
    t(G)\ge \frac{\lambda_1(G)(\lambda_1^2(G) - s^-(G))}{6}. 
\end{equation}

One can compare  \eqref{eq:triangle_s_minus_rearranged} with the following well-known triangle counting result of Bollob\'{a}s and Nikiforov \cite{bollobas2007cliques}: for any graph $G$,  
\begin{equation}\label{eq:triangle_BN}
    t(G)\ge \frac{\lambda_1(G)(\lambda_1^2(G) - m)}{3}. 
\end{equation}
Inequality \eqref{eq:triangle_BN} was also independently observed by Cioab\u{a}, Feng, Tait and Zhang \cite{Cioaba_Feng_Tait_Zhang_2020}, and the equality case was analyzed by Ning and Zhai \cite{Ning_Zhai_2023}; also refer \cite{Li_Liu_Zhang_2025_Nosal}. 

Owing to the difficulty of proving Conjecture \ref{conj:splus_omega} and \ref{conj:vertex_local_wilf}, it is of interest to determine a non-trivial function $f: \mathbb{N} \rightarrow \mathbb{N}$ such that 
\[\sqrt{s^+(G)}\le \sum_{v\in V(G)}\left(1 - \frac{1}{f(c(v))}\right).\]
We show that $f(\omega) = 2\omega$ is sufficient.

\begin{proposition}
For a graph $G$, we have 
    \[ \sqrt{s^+(G)} \le \sum_{v\in V(G)}\left(1-\frac{1}{2c(v)}\right). \]
\end{proposition}

\begin{proof}
    Let $y\in \mathbb{R}^{\vert V(G) \vert}$ be the vector such that $y_v = \sqrt{\frac{c(v)-1}{c(v)}}$ for all $v \in V(G)$. Now 
\begin{align*}
     F_G(y) &=  \sum_{uv \in E(G)} 2\sqrt{\frac{c(u)c(v)}{(c(u) -1)(c(v)-1)}}y_uy_v = \sum_{uv \in E(G)} 2 = 2m.
\end{align*}
Using Lemma \ref{lemma:weighted_M_S_2}, we have 
\[F_G(y) \leq \left(\sum_{v\in V(G)} y_v\right)^2 = \left(\sum_{v \in V(G)} \sqrt{1-\frac{1}{c(v)}} \right)^2.\]
Since $s^+ + s^- = 2m$, we have 
\[\sqrt{s^+} \leq \sqrt{2m} \leq \left( \sum_{v\in V(G)}\sqrt{1-\frac{1}{c(v)}} \ \right) \leq  \sum_{v\in V(G)}\left(1-\frac{1}{2c(v)}\right). \qedhere \] 
\end{proof}


For Conjecture \ref{conj:splus_omega}, we can improve the function to $f(\omega) = 2\omega - 2$, thus proving Theorem \ref{thm:splus_weak_2omegaminus2}.

\begin{proof}[Proof of Theorem \ref{thm:splus_weak_2omegaminus2}]
For every vertex $v \in V(G)$, the subgraph induced by the open neighbourhood of $v$, i.e., $G[N(v)]$ is $K_{c(v)}$-free and has at most $\frac{\deg(v)^2}{2}\left(1 - \frac{1}{c(v)-1}\right)$ edges by Tur\'{a}n's Theorem. Using Lemma \ref{lemma:triangle_counting}, we get 
\begin{equation}\label{eq:triangle_1}
  t(G)\le \sum_{v \in V}\frac{\deg(v)^2}{6}\left(1 - \frac{1}{c(v)-1}\right).  
\end{equation}
Since $\sum_{v\in V(G)}\deg(v)^2 \le n\lambda_1^2$ for any graph $G$, we get 
\[t(G)\le  \sum_{v \in V}\frac{\deg(v)^2}{6}\left(1 - \frac{1}{\omega -1}\right)\le \frac{n\lambda_1^2}{6}\left(1 - \frac{1}{\omega-1}\right).\]
Using Theorem \ref{thm:splus_triangles}, we see
\begin{align*}
    \sqrt{s^+} \leq \frac{n}{2} + \frac{3t(G)}{\lambda_1^2} \le n\left(1 - \frac{1}{2\omega-2}.\right)
\end{align*}
Next, we improve the above bound even further. Rearranging inequality \ref{eq:triangle_s_minus} and noting that $s^- = 2m - s^+$, we get
\begin{align*}
    s^+ \leq 2m - \lambda_1^2 + \frac{6t(G)}{\lambda_1}.
\end{align*}
Since $2m\leq n\lambda_1$ and substituting the bound for $t(G)$, we have
\[s^+ \leq n\lambda_1\left(2 - \frac{1}{\omega-1}\right) - \lambda_1^2\]
Wilf's inequality gives us $\lambda_1 \leq n\left(1 - \frac{1}{\omega}\right)$, so maximing the function on the right gives
\begin{align*}
    s^+ \le n^2 \left(1 - \frac{1}{\omega} - \frac{1}{\omega^2}\right). 
\end{align*}
This completes the proof. 
\end{proof}

For regular graphs, the following is immediate from the proof of Theorem \ref{thm:splus_weak_2omegaminus2}. 

\begin{proposition}
If $G$ is a regular graph, then 
 \[ \sqrt{s^+(G)} \le \sum_{v\in V(G)}\left(1-\frac{1}{2c(v)-2}\right). \]
\end{proposition}

\subsection{Diamond-free graphs}

Conjecture \ref{conj:vertex_local_wilf} is equivalent to Conjecture \ref{conj:splus_omega} for triangle-free graphs and hence true. We make further progress by proving Conjecture \ref{conj:vertex_local_wilf} for diamond-free graphs.

\begin{theorem}\label{thm:local_wilf_diamond_free}
Let $G$ be a diamond-free graph of order $n\ge 42$. Then
\[\sqrt{s^+(G)} \leq \sum_{v \in V(G)} \left(1 - \frac{1}{c(v)}\right).\] 
\end{theorem}

\begin{proof}
Observe that we can write the right-hand side of the inequality as follows: 
\begin{equation}\label{eq:diamond_free_1}
  \sum_{v \in V(G)} \left(1 - \frac{1}{c(v)}\right) = \frac{n}{2} + \sum_{v \in V(G)} \left(\frac{1}{2} - \frac{1}{c(v)}\right).  
\end{equation}
Since $c(v) \ge 2$ for every vertex $v$, this gives that the right-hand side of the inequality in the assertion is always at least $n/2$. 

If $G$ is triangle-free, then $t(G)=0$ and thus Theorem \ref{thm:splus_triangles} implies $\sqrt{s^+} \leq \frac{n}{2}$. So, assume $t(G)\ge 1$. Since $G$ is diamond-free, it follows that $\omega(G)=3$ and $c(v) \le 3$ for all vertices $v \in V$. 

If $\lambda_1(G)\le \frac{n}{4}$, then $\sqrt{s^+} \le \sqrt{n\lambda_1} \le \frac{n}{2}$, and we are done. So assume $\lambda_1\ge \frac{n}{4}$.

Since $G$ is diamond-free, for every vertex $v \in V(G)$, the induced subgraph $G[N(v)]$ is $P_3$-free and has at most $\frac{\deg(v)}{2}$ edges by Theorem \ref{thm:erdos_gallai}. Using Lemma \ref{lemma:triangle_counting}, we get
\[t(G)\le \sum_{v\in V(G)}\frac{\deg(v)}{2} = \frac{m}{3}.\] 
Applying Theorem \ref{thm:splus_triangles} and using the fact that $2m\le n\lambda_1$, we have
\begin{equation}\label{eq:diamond_free_2}
  \sqrt{s^+} \leq  \frac{n}{2} + \frac{3t(G)}{\lambda^2_1} \leq  \frac{n}{2} + \frac{m}{\lambda^2_1} \le \frac{n}{2} + 2.  
\end{equation}

Define the set of \emph{triangular vertices} in $G$ by 
\[ TV(G) =\{v\in V(G): c(v) = 3\},\]
and let $tv(G) = |TV(G)|$. We consider the following cases:

\textbf{Case 1:} $tv(G)\ge 12$. 

Then 
\[\sum_{v \in V(G)} \left(1 - \frac{1}{c(v)}\right) \ge \frac{n}{2} + 2,\]
and the assertion holds by \eqref{eq:diamond_free_2}.

\textbf{Case 2:} $tv(G)\le 11$.

Note that $3t(G) \le {\frac{tv(G)(tv(G)-1)}{2}}\le \frac{11\cdot 10}{2}$, implying $t(G)\le \lfloor \frac{110}{6} \rfloor = 18$. Using \eqref{eq:diamond_free_2}, we have 
\[ \sqrt{s^+} \leq  \frac{n}{2} + \frac{3t(G)}{\lambda^2_1} \leq  \frac{n}{2} + \frac{3\cdot 18}{n^2/4^2} \le \frac{n}{2} + \frac{1}{2},\]
whenever $n\ge 42$. Since $t(G)\ge 1$, by \eqref{eq:diamond_free_1} we have 
\[\sum_{v \in V(G)} \left(1 - \frac{1}{c(v)}\right) \ge \frac{n}{2} + \frac{1}{2},\]
completing the proof.
\end{proof}

\subsection{$F$-free graphs}

\begin{theorem}\label{thm:ESS_general_F_free_chi_3}
For any $F$-free graph $G$ with chromatic number $\chi(F) = 3$, we have
    \[\sqrt{s^+(G)} \le n \left(\frac{1}{2} + o(1)\right).\]
\end{theorem}

\begin{proof}
If $\lambda_1\le \frac{n}{4}$, then $\sqrt{s^+}\le \sqrt{n\lambda_1}\le \frac{n}{2}$, and the assertion holds. So assume $\lambda_1 \ge \frac{n}{4}$. 

Using Corollary \ref{cor:triangles_F_free}, we have $t(G) \le o(n^3)$. Using Theorem \ref{thm:splus_triangles}, we have 
\[        \sqrt{s^+} \leq \frac{n}{2} + \frac{3t(G)}{\lambda^2_1} \le \frac{n}{2} + \frac{o(n^3)}{n^2} = n\left(\frac{1}{2} + o(1) \right).\qedhere\]
\end{proof}

\subsection{Random graphs}

Here, we prove Conjecture \ref{conj:vertex_local_wilf} for the random graph $G(n, \frac{1}{2})$. 

\begin{theorem}\label{thm:local_wilf_random}
    For the random graph $G = G(n, \frac{1}{2})$, we have 
    \[ \sqrt{s^+(G)} \le \sum_{v\in V(G)} \frac{c(v)-1}{c(v)}. \]
\end{theorem}

\begin{proof}
 Note that for any vertex $v\in V(G)$, $\deg(v) = \left(\frac{1}{2} + o(1)\right)n$ a.a.s, and so 
    \[c(v) = \omega(G[N(v)]) = \omega(G(n/2, 1/2)) = \frac{2\log(n/2)}{\log(2)}\]
    a.a.s. Using Theorem \ref{thm:random_graph_square_energy}, we see that a.a.s  
\[ \sqrt{s^+} \le  n\sqrt{\frac{3}{8} + o(1)} \le n \left(1 - \frac{\log(2)}{2\log(n/2)}\right) = \sum_{v\in V(G)} \frac{c(v)-1}{c(v)}.\qedhere\]
\end{proof}

\section{Generalizations of Bollob\'{a}s-Nikiforov conjecture}
\label{sec:BN}

\subsection{Diamond-free graphs}

In this section, we prove Conjecture \ref{conj:local_BN} for diamond-free graphs, which is the next non-trivial case after triangle-free graphs. Kumar and Pragada \cite{Bollobas_Nikiforov_few_triangles_2025} have shown that if diamond-free graphs with at least $34$ edges, then Conjecture \ref{conj:BN} holds. We first show here that Conjecture \ref{conj:BN} holds for all diamond-free graphs, as this is required later in the proof.

\begin{theorem}\label{thm:BN_diamond_free}
    Let $G$ be a diamond-free graph. Then 
    \[\lambda_1^2(G) + \lambda^2_2(G) \leq 2m\left(1 - \frac{1}{\omega} \right).\]
\end{theorem}
\begin{proof}
We can assume that $ n\ge 5$. Since $G$ is not a complete graph, $\lambda_2\ge 0$.
 Clearly, $\omega \le 3$. The assertion is known to be true for triangle-free graphs, so we can assume that $\omega = 3$. Suppose to the contrary that $\lambda_1^2+ \lambda_2^2 > \frac{4m}{3}$. Then,
\[\sum_{\lambda_i<0}\lambda_i^2 \le 2m - (\lambda_1^2 + \lambda_2^2) < \frac{2m}{3}.\]
Now, we count the number of triangles in $G$ in two different ways. As in the proof of Theorem \ref{thm:local_wilf_diamond_free}, we have $t(G)\le \frac{m}{3}$. Thus,
    \begin{align*}
      2m \ge 6t(G) &= \sum_{i=1}^n\lambda_i^3 \\
      &\ge \lambda_1^3 + \lambda_2^3 - \sum_{\lambda_i<0} |\lambda_i^3| \\
      &\ge \frac{(\lambda_1^2 + \lambda_2^2)^{3/2}}{\sqrt{2}} - \left(\sum_{\lambda_i<0} |\lambda_i^2|\right)^{3/2} \\
      &> \frac{1}{\sqrt{2}}\left(\frac{4m}{3}\right)^{3/2} - \left(\frac{2m}{3}\right)^{3/2} 
    \end{align*}
This is a contradiction whenever $m\ge 14$. If $m \le 13$, then $n \le 14$ since $G$ is connected, and the assertion is verified using a computer \href{https://github.com/Shivaramkratos/Code-for-Localization/blob/main/diamond_free_check.sage}{here}. 
\end{proof}

In \cite{Bollobas_Nikiforov_few_triangles_2025}, the following general bound was established for $\lambda_1^2 + \lambda_2^2$ in terms of the number of triangles.

\begin{theorem}[\cite{Bollobas_Nikiforov_few_triangles_2025}]\label{thm:BN_triangles}
    Let $G$ be a graph of size $m$ and number of triangles $t(G)$. Then
    \[\lambda_1^2(G) + \lambda_2^2(G) < m + \left(3t(G)\right)^{2/3}.\]
\end{theorem}

For diamond-free graphs, we improve the coefficient of $t(G)$ in the upper bound given in Theorem \ref{thm:BN_triangles}. Our proof relies on majorization, refer to Subsection \ref{subsec:majorization}.

\begin{lemma}\label{lemma:BN_triangles_improved_diamond_free}
    Let $G$ be a diamond-free graph of size $m$ and number of triangles $t(G)$. Then 
    \[\lambda_1^2(G) + \lambda_2^2(G) \leq m + \left(\frac{3}{\sqrt{2}}t(G)\right)^{2/3}.\]
\end{lemma}

\begin{proof}
We can assume that $n\ge 5$ and $\lambda_2\ge 0$. If $\lambda_1^2 + \lambda_2^2 \leq m$, then we are done. So suppose $\lambda_1^2 + \lambda_2^2 = m + \delta$, for some $\delta>0$. By Theorem \ref{thm:BN_diamond_free}, $\delta \leq m/3$, implying 
\begin{equation}\label{eq:diamond_triangle_improved_1}
    \lambda_1^2 \ge \frac{\lambda_1^2 + \lambda_2^2}{2}  = \frac{m + \delta}{2} \ge 2\delta.
\end{equation}
Since $\sum_{i=1}^n \lambda_i^2 = 2m$, we have
\[ 2(\lambda_1^2 + \lambda_2^2) =  2m + 2\delta = \sum_{i=1}^n \lambda_i^2 + 2\delta,\]
which implies
\begin{equation}\label{eq:diamond_triangle_improved_2}
\lambda_1^2 + \lambda_2^2 = \sum_{i=3}^n \lambda_i^2 + 2\delta \geq \sum_{\lambda_i < 0}^n\lambda_i^2 + 2\delta.
\end{equation}
Let $x = (\lambda_1^2, \lambda_2^2,0,\dots,0)^T$ and $y = (\lambda_n^2, \lambda_{n-1}^2,\dots, \lambda^2_{n-n^{-}+1},2\delta)^T$ in $\mathbb{R}^{n^{-}+1}$, where $n^-$ denotes the number of negative eigenvalues of $G$. Clearly, $\lambda_1^2\ge \max\{\lambda_n^2, 2\delta\}$ by \eqref{eq:diamond_triangle_improved_1}. Then, using \eqref{eq:diamond_triangle_improved_2} we see that $y \prec_w x$. Applying Theorem \ref{thm:majorization} with $p=\frac{3}{2}$, we get $\Vert x \Vert^{3/2}_{3/2} \geq \Vert y \Vert^{3/2}_{3/2}$, that is,
\[\lambda_1^3 +  \lambda_2^3 \ge \sum_{\lambda_i<0}^n|\lambda_i|^3 + (2\delta)^{3/2}.\]
This implies that
\begin{align*}
    6t(G) &= \sum_{i=1}^n \lambda_i^3 \geq \lambda_1^3 + \lambda_2^3 - \sum_{\lambda_i<0}^n|\lambda_i|^3\\ &\ge (2\delta)^{3/2} = 2\sqrt{2}\big(\lambda_1^2 + \lambda_2^2- m\big)^{3/2}.
\end{align*}
Upon rearrangement, we get the desired inequality.
\end{proof}

\begin{theorem}\label{thm:local_BN_diamond_free}
    Let $G$ be a diamond-free graph of size $m$ and number of triangles $t(G)$. If $t(G)\notin \{1,2,3,4\}$, then
    \[\lambda_1^2 + \lambda_2^2 \leq \sum_{e\in E(G)} 2\left(1-\frac{1}{c(e)}\right).\]
\end{theorem}

\begin{proof} If $t(G)=0$, then the assertion is known to be true. So assume $t(G)\ge 5$. Since $G$ is diamond-free, for any $e \in E(G)$,  $c(e) \in \{2, 3\}$. Moreover, if $c(e) = 3$, then $e$ is contained in a unique triangle. Thus, we have
\begin{equation}
   \sum_{e\in E(G)} 2\left(1-\frac{1}{c(e)}\right) = m + \sum_{e\in E(G)} \left(1-\frac{2}{c(e)}\right) = m + t(G).
\end{equation}
Since $t(G) \geq 5$ and using Lemma \ref{lemma:BN_triangles_improved_diamond_free}, we have 
    \[\lambda_1^2 + \lambda_2^2 \leq m + \left(\frac{3}{\sqrt{2}}t(G)\right)^{2/3} \leq m + t(G). \qedhere\]
\end{proof}

When $t(G)\in \{1,2,3,4\}$, Lemma \ref{lemma:BN_triangles_improved_diamond_free} gives an upper bound which is only slightly worse than $m+t(G)$.

\subsection{$F$-free graphs}

\begin{theorem}\label{thm:BN_F_free_chi_3}
Let $F$ be a fixed graph with $\chi(F)=3$. Then every $F$-free graph $G$ of size $m$ satisfies
\[\lambda_1^2(G) + \lambda_2^2(G) \le m(1+o(1)).\]
\end{theorem}

\begin{proof} 
Using Corollary \ref{cor:triangles_F_free}, we have $t(G)\le o(m^{3/2})$. Using Theorem \ref{thm:BN_triangles}, we get 
\[\lambda_1^2(G) + \lambda_2^2(G) \le m + (3t(G))^{2/3} \le m + o(m) = m(1 + o(1)).\qedhere\]
\end{proof}

\subsection{Random graphs}

Here, we verify Conjecture \ref{conj:local_BN} for the random graph $G(n, \frac{1}{2})$, thus showing that it is true for almost all graphs.  

\begin{theorem}\label{thm:local_BN_random}
For the random graph $G = G(n, \frac{1}{2})$, we have 
    \[\lambda_1^2(G) + \lambda_2^2(G) \le \sum_{e\in E(G)} 2\left(1-\frac{1}{c(e)}\right).\]
\end{theorem}

\begin{proof}
 Observe that for any edge $uv\in E(G)$, we have 
 \[c(uv) = \omega(G[N(u)\cap N(v)]) = \omega(G(n/4, 1/2)) = \frac{2\log (n/4)}{\log 2}\]
 a.a.s. Using Theorem \ref{thm:random_graph_square_energy}, we have a.a.s.
 \[ s^+ \le  n^2\left(\frac{3}{8} + o(1)\right) \leq \frac{n^2}{2}\left( 1 - \frac{\log 2}{2\log (n/4)}\right) = \sum_{e\in E(G)} 2\left(1-\frac{1}{c(e)}\right).  \qedhere\]
\end{proof}

\section{Walk inequalities} 
\label{sec:Walk_inequalities}

In this section, we prove a localized version of Nikiforov's Walk Inequality. First, we observe the following recursive inequality.

\begin{lemma}\label{lemma:rescursion_walks} For any integer $r\ge 1$, we have 
\[ w_{2r}(G) \le \left( \sum_{v\in V(G)}w_r(v)\sqrt{\frac{c(v)-1}{c(v)}}\right)^2.\]
\end{lemma}

\begin{proof}
Let $y\in \mathbb{R}^{V(G)}$ be the vector such that $y_v = w_r(v)\sqrt{\frac{c(v)-1}{c(v)}}$ for all $v \in V(G)$. Now 
\begin{align*}
     F_G(y) &=  \sum_{uv \in E(G)} 2\sqrt{\frac{c(u)c(v)}{(c(u) -1)(c(v)-1)}}y_uy_v = \sum_{uv \in E(G)} 2{w_r(u)w_r(v)} = w_{2r}(G).
\end{align*}
Using Lemma \ref{lemma:weighted_M_S_2}, we have $F_G(y) \leq \left(\sum_{v\in V(G)} y_v\right)^2$ which proves the assertion.
\end{proof}

We are now ready to prove Theorem \ref{thm:partial_walks}.

\begin{proof}[Proof of Theorem \ref{thm:partial_walks}]
Using Nikiforov's Walk Inequality (Theorem \ref{thm:Nikiforov_walk}) and Lemma \ref{lemma:rescursion_walks}, we have 
\begin{equation*}
\lambda_1^{2r} \le \frac{\omega-1}{\omega}\,w_{2r}(G) \le  \frac{\omega-1}{\omega}\left(\sum_{v\in V(G)} w_r(v)\sqrt{\frac{c(v)-1}{c(v)}}\right)^2
\end{equation*}
Taking square roots gives the desired result. Equality case is easily seen to be the same as in Theorem \ref{thm:Nikiforov_walk}.
\end{proof}

\section*{Acknowledgement}

M. Rajesh Kannan acknowledges financial support from the ANRF-CRG India and SRC, IIT Hyderabad. The authors thank Clive Elphick for several helpful comments. We also thank Rajat Adak, Yongtao Li and Xiao-Dong Zhang for suggesting some references. 

\bibliographystyle{plainurl}
\bibliography{references.bib}

\vspace{0.4cm}
\noindent M. Rajesh Kannan, Email: {\tt rajeshkannan@math.iith.ac.in, rajeshkannan1.m@gmail.com}\\
\textsc{Department of Mathematics, Indian Institute of Technology Hyderabad, Sangareddy 502285, India}\\[1pt]

\noindent Hitesh Kumar, Email: {\tt hitesh.kumar.math@gmail.com}, {\tt hitesh\_kumar@sfu.ca}\\
\textsc{Department of Mathematics, Simon Fraser University, Burnaby, Canada}\\[1pt]

\noindent Shivaramakrishna Pragada, Email: {\tt shivaramakrishna\_pragada@sfu.ca, shivaramkratos@gmail.com}\\
\textsc{Department of Mathematics, Simon Fraser University, Burnaby, Canada}

\end{document}